\documentclass{amsart}

\usepackage{amsthm,amsmath,amssymb,url}
\newtheorem{thm}{Theorem}[section]
\newtheorem{lem}[thm]{Lemma}
\newtheorem{cor}[thm]{Corollary}
\theoremstyle{remark}
\newtheorem*{rmk}{Remark}

\newtheorem*{eg*}{Example}

\newcommand{\coeff}[1]{\langle #1 \rangle}

\begin{document}

\title[$q$-Faulhaber coefficients]{A combinatorial interpretation of Guo and
  Zeng's $q$-Faulhaber coefficients}

 \thanks{Institut f\"ur Statistik und Decision Support\\
   \url{martin.rubey@univie.ac.at}\\
   \url{http://www.mat.univie.ac.at/~rubey}}

\begin{abstract}
  Recently, Guo and Zeng discovered $q$-analogues of Faulhaber's formulas for
  the sums of powers. They left it as an open problem to extend the
  combinatorial interpretation of Faulhaber's formulas as given by Gessel and
  Viennot to the $q$ case. In this note we will provide such an interpretation.
\end{abstract}
\maketitle

\section{Introduction}
In the early seventeenth century, Johann Faulhaber considered the sums of
powers $S_{m,n}=\sum_{k=1}^n k^m$ and provided formulas for the coefficients
$f_{m,k}$ in
\begin{equation*}
  S_{2m+1,n}=\sum_{k=1}^m f_{m,k}\left(\frac{n(n+1)}{2}\right)^{k+1}.
\end{equation*}

A combinatorial interpretation of these coefficients was given by Gessel and
Viennot\cite{GesselViennot1989}. Recently, Guo and Zeng\cite{GuoZeng2005},
continuing work of Schlosser\cite{Schlosser}, Warnaar\cite{Warnaar} and Garrett
and Hummel\cite{GarrettHummel} were able to find $q$-analogues of Faulhaber's
formulas.

More precisely, setting $[k]=\frac{1-q^k}{1-q}$ and
$[k]!=\prod_{i=1}^k\frac{1-q^i}{1-q}$, they proved the following result:
\begin{thm}
  Let $m,n\in\mathbb{N}$ and set
  \begin{equation}\label{eq:1}
    S_{m,n}(q)=\sum_{k=1}^{n}\frac{[2k]}{[2]} [k]^{m-1}q^{\frac{m+1}{2}(n-k)}.
  \end{equation}
  Then there exist polynomials $P_{m,k}\in\mathbb{Z}[q]$ such that
  \begin{equation}\label{eq:2}
    S_{2m+1,n}(q) =\sum_{k=0}^m (-1)^{m-k} \frac{[k]!}{[m+1]!}P_{m,m-k}
                                q^{(m-k)n} \frac{([n][n+1])^{k+1}}{[2]}.
  \end{equation}
\end{thm}

They also gave an explicit formula for the polynomials $P{m,k}(q)$, and asked
for a combinatorial interpretation of the coefficients of these polynomials. It
is the purpose of this paper to answer this question. 

\section{A combinatorial interpretation of the $q$-Faulhaber coefficients}

In this Section we will exhibit a surprisingly simple combinatorial
interpretation of the polynomials $P_{m,k}$ which also yields some interesting
properties as easy consequences. We follow exactly the arguments given by
Gessel and Viennot\cite[Theorem~29]{GesselViennot1989}:

\begin{thm}\label{thm:inverse-matrix}
  Let $h_{2m-k}(\{1,q\}^{k-m+1})$ denote the $(2m-k)$\textsuperscript{th}
  complete homogeneous function in $2(k-m+1)$ variables, half of which are
  specialised to~$q$, the others to~$1$. Thus, for $m\leq k$ we have
  $h_{2m-k}(\{1,q\}^{k-m+1})
  =\sum_{i=0}^{2m-k}\binom{k-m+i}{k-m}\binom{m-i}{k-m}q^i$, for $m>k$ it
  vanishes.
  
  Then the inverse of the matrix
  \begin{equation*}
    \left(h_{2m-k}(\{1,q\}^{k-m+1})\right)_{k,m\in\{0\dots n\}}
  \end{equation*}
  is the matrix
  \begin{equation*}
    \left((-1)^{k-m}\frac{[m]!}{[k+1]!}P_{k,k-m}\right)_{k,m\in\{0\dots n\}}.
  \end{equation*}
  Note that both matrices are lower triangular.
\end{thm}

The proof rests on the following Lemma, which might be interesting per se. When
$q=1$ it reduces to a simple application of the binomial theorem. In the
general case however, it turns out to be considerably more difficult.

\begin{lem}\label{lem:1}
  For any $m$ and $l$ we have
  \begin{multline}\label{eq:3}
    \left(\frac{[l][l+1]}{q^l}\right)^{m+1}
   -\left(\frac{[l-1][l]}{q^{l-1}}\right)^{m+1}\\
   =[2]\sum_k h_{m-2k}(\{1,q\}^{k+1})
              \frac{[2l]}{[2]}[l]^{2(m-k)}
              q^{-l(m-k+1)}.
  \end{multline}
\end{lem}
\begin{proof}
  Multiplying both sides with $\left(\frac{q^l(1-q)}{[l]}\right)^{m+1}$,
  and replacing $q^l$ with $x$, we arrive at the following equivalent identity:
  \begin{multline}
    \label{eq:5}
    (1-qx)^{m+1}-(q-x)^{m+1}\\
    =\sum_{k} h_{m-2k}(\{1,q\}^{k+1})(1-q)^{2k+1}(1-x)^{m-2k}(1+x)x^k.
  \end{multline}
  We will show that the coefficients of $q^r x^s$ are the same on both sides.
  To this end, we first extract the coefficient of $q^r$ in
  $h_{m-2k}(\{1,q\}^{k+1})(1-q)^{2k+1}$ which for $k\leq \lfloor\frac m
  2\rfloor$ turns out to be
  \begin{equation*}
    \sum_{i=0}^{m-2k}(-1)^{r-i}\binom{k+i}{k}\binom{m-k-i}{k}\binom{2k+1}{r-i}.
  \end{equation*}
  Transforming this sum into hypergeometric notation -- using for example
  Christian Krattenthaler's package {\tt hyp.m}\cite{HYPHYPQ} -- we find that it
  equals
  \begin{equation*}
    (-1)^r \binom{m-k}{k} \binom{2k+1}{r}\,
    {}_3 F_2 \!\left [ 
    \begin{matrix} 
       { -r, 1 + k, 2k - m}\\ 
       { k-m, 2 + 2k - r}
    \end{matrix}; 1\right ],
  \end{equation*}
  which is summable by \cite[(2.3.1.3); Appendix (\rm III.2)]{Slater}. After
  some simplification we obtain
  \begin{equation*}
    (-1)^{r+m} \binom{m+1}{r} \binom{r-k-1}{m-2k}
  \end{equation*}
  as closed form for the sum. Since the sign and the first binomial coefficient
  does not involve the summation index $k$, it remains to evaluate the
  coefficient of $x^s$ in
  \begin{equation*}
    \sum_{k=0}^{\lfloor\frac{m}{2}\rfloor}\binom{r-k-1}{m-2k}(1-x)^{m-2k}(1+x)x^k.
  \end{equation*}
  Obviously it is sufficient to find a closed form for
  \begin{equation*}
    \coeff{x^s} \sum_{k=0}^{\lfloor\frac{m}{2}\rfloor}
                    \binom{r-k-1}{m-2k}(1-x)^{m-2k}x^k
    = \sum_{k=0}^{\lfloor\frac{m}{2}\rfloor}
                 (-1)^{s-k}\binom{r-k-1}{m-2k}\binom{m-2k}{s-k}.
  \end{equation*}
  Now we have to distinguish several cases. For $2k \leq m$ we observe that the
  product of binomial coefficients $\binom{r-k-1}{m-2k}\binom{m-2k}{s-k}$ does
  not vanish only if
  \begin{equation*}
    \min(r,m-r+1)\leq k \leq \min(s, m-s).   
  \end{equation*}
  Therefore, the sum above must be zero for 
  \begin{gather*}
    s<r\text{ or } s>m-r\text{, if } r\leq\lfloor\frac{m}{2}\rfloor\\
  \intertext{and}
    s\geq r\text{ or } s\leq m-r\text{, if } r>\lfloor\frac{m}{2}\rfloor.
  \end{gather*}
  If the pair $r$ and $s$ happens to be outside of this region, the bounds of
  summation are natural, so we can write the sum as a hypergeometric series.
  This time, it equals
  \begin{equation*}
    (-1)^s \binom{r-1}{m} \binom{m}{s}\,
    {}_3 F_2 \!\left [
    \begin{matrix} 
       { -s, -m + s, 1}\\ 
       { 1 - r, -m + r}
    \end{matrix} ; 1\right ],
  \end{equation*}
  which is not directly summable. However, the terminating form of the
  transformation in \cite[Ex. 7, p. 98]{Bailey} is applicable. Unfortunately,
  the parameter that causes the sum to terminate cancels, so we have to
  evaluate
  \begin{equation*}
    {}_3 F_2 \!\left [
    \begin{matrix} 
       { -s, -m + s, 1+\varepsilon}\\ 
       { 1 - r, -m + r}
    \end{matrix} ; 1\right ]
  \end{equation*}
  instead. Again, we have to distinguish between $s\geq r$ and $s>m-r$,
  however, the only change in the computation amounts to exchanging the two
  lower parameters of the hypergeometric expression above.  Applying the above
  mentioned transformation twice and then sending $\varepsilon$ to zero, we
  arrive at an expression involving
  \begin{gather*}
    {}_2 F_1 \!\left [ 
    \begin{matrix} 
       { r - s, r - m + s}\\ 
       { r - m}
    \end{matrix} ; {\displaystyle 1}\right ],
  \intertext{or, for $s>m-r$ }
    {}_2 F_1 \!\left [ 
    \begin{matrix} 
       { 1 - r +m - s, 1 - r + s}\\ 
       { 1 - r}
    \end{matrix} ; {\displaystyle 1}\right ]
  \end{gather*}
  to which we can apply \cite[Appendix (III.4)]{Slater}. After some
  simplification we see that the sum evaluates to
  \begin{align*}
    (-1)^{m+r+s}  &\text{ if } r\leq s\leq m-r \text{ and } 
                               r\leq\lfloor\frac{m}{2}\rfloor\\
  \intertext{and}
    (-1)^{m+r+s+1}&\text{ if } m-r < s < r \text{ and } 
                               r>\lfloor\frac{m}{2}\rfloor.
  \end{align*}

  Putting all the pieces together we obtain
  \begin{align*}
    (-1)^r\binom{m+1}{r} &\text{ if } s=r\neq m-r+1,\\
    (-1)^{r+m}\binom{m+1}{r}&\text{ if } s=m-r+1\neq r\text{ and}\\
    0&\text{ otherwise}
  \end{align*}
  for the coefficient of $q^r x^s$ in the right hand side of \eqref{eq:3}. This
  is easily seen to be equal to the coefficient of $q^r x^s$ in the left hand
  side of \eqref{eq:3}.
\end{proof}
\begin{rmk}
  It seems that an alternative proof could be given based on the
  Wilf-Zeilberger method. In any case, it would be interesting to see the given
  identity as a special case of a more general one. In particular, is there a
  $q$-binomial theorem that could be applied to $[l+1]^{m+1}$?

  A second question poses itself by looking at the corresponding section in
  the paper of Gessel and Viennot: In the standard case, a very similar
  identity can be derived for the Sali\'e numbers, i.e., the coefficients of
  $\left(n(n+1)\right)^k$ in $\sum_{k=1}^n (-1)^{n-k} k^m$. However, extending
  the ideas of this paper in a straightforward fashion does not seem to work.
  
  Schlosser\cite{Schlosser} considered various $q$ analogues of these
  alternating sums and derived formulas for $m\leq 4$, the most plausible being
  \begin{gather*}
    T_{m,n}(q)=\sum_{k=1}^{n}\frac{[2k]}{[2]}[k]^{m-1}(-q^{\frac{m+1}{2}})^{n-k}
  \end{gather*}

  Standard computer algebra packages are able to find formulas for greater
  values of $m$. It turns out that the coefficient of $\left([n][n+1]\right)^k$
  in $T_{2m,n}$ is of the form 
  \begin{equation*}
    \frac{(1+q^{n+\frac 1 2})(1+q^{\frac 1 2})^{m-k}}
         {(1+q^{\frac{2m+1} 2})\dots(1+q^{\frac{2k+1} 2})} 
    (-q^n)^{m-k} 
    g_{k,m}(q),
  \end{equation*}
  where $g_{k,m}(q)$ is a polynomial in $q$. For small $m$ these polynomials
  are listed in Table~\ref{tab:Salie}. Unfortunately, they do not have
  nonnegative coefficients, which indicates that this is not the \lq right\rq\ 
  $q$-analogue of the Sali\'e numbers. It would be particularly nice to have a
  common $q$-analogue of Faulhaber and Sali\'e numbers. Note that Guo and
  Zeng\cite{GuoZeng2005} proposed a more general $q$ analogue of Faulhaber's
  numbers, but these also fail to have nonnegative coefficients.

  \begin{table}[h]
    \begin{equation*}
      \begin{array}{|l|r|r|r|r|r|r|}
        \hline
        k\setminus m & 1 & 2 & 3                 & 4\\\hline
        1            & 1 & 1 & 2q-q^\frac 1 2 +2 & 
(5q^2-q^\frac 3 2 +9q-q^\frac 1 2+5)(q-q^\frac 1 2 + 1)\\\hline
        2            &   & 1 & 2q-q^\frac 1 2 +2 & 
(5q^2-q^\frac 3 2 +9q-q^\frac 1 2+5)(q-q^\frac 1 2 + 1)\\\hline
        3            &   &   & 1                 & 
3q^2-2q^\frac 3 2 +4q-2q^\frac 1 2 + 3\\\hline
        4            &   &   &                   & 1\\\hline
      \end{array}
    \end{equation*}
    \caption{$g_{k,m}(q)$ for $m\leq 3$\label{tab:Salie}}
  \end{table}

  Finally, although the appearance of the complete homogeneous symmetric
  functions is natural, the specialisation involved seems to be interesting and
  might deserve more attention.
\end{rmk}

\begin{proof}[Proof of Theorem~\ref{thm:inverse-matrix}]
  Summing Equation~\eqref{eq:3} on $l$ from $1$ to $n$, observing that its left
  hand telescopes and using Equation~\eqref{eq:1} on its right hand side we
  obtain
  \begin{equation*}
    \begin{split}
      \left(\frac{[n][n+1]}{q^n}\right)^{m+1}
    &=[2]\sum_{k} h_{m-2k}(\{1,q\}^{k+1}) S_{2(m-k)+1,n}q^{-n(m-k+1)}.
    \end{split}
  \end{equation*}
  Plugging in Equation~\eqref{eq:2} and exchanging the order of summation, the
  right hand side becomes
  \begin{equation*}
    \sum_{l}\sum_{k\geq l} h_{m-2k}(\{1,q\}^{k+1}) 
                          (-1)^{m-k-l}\frac{[l]!}{[m-k+1]!}P_{m-k,m-k-l}
                          \left(\frac{[n][n+1]}{q^n}\right)^{l+1}.
  \end{equation*}
  Comparing coefficients of $\left(\frac{[n][n+1]}{q^n}\right)^{l+1}$ we see
  that the two matrices in question are indeed inverses.
\end{proof}

We also copy a simple lemma from Gessel and Viennot\cite{GesselViennot1989},
that follows easily from the formula for the entries of the inverse of a
matrix:
\begin{lem}
  Let $\left(A_{i,j}\right)_{i,j\in\{1,2,\dots,m\}}$ be an invertible lower
  triangular matrix and let $B$ be its inverse. Then for $0\leq k\leq n\leq m$
  we have
  \begin{equation*}
    B_{n,k}=\frac{(-1)^{n-k}}{A_{k,k},A_{k+1,k+1},\dots,A_{n,n}}
            \det\left(A_{k+i+1,k+j}\right)_{i,j\in\{0,1,\dots,n-k-1\}}.
  \end{equation*}
\end{lem}

Finally we can announce our main theorem:
\begin{thm}
  \begin{equation}\label{eq:4}
     P_{m,k}=\det\left(h_{m-k-i+2j-1}(\{1,q\}^{i-j+2})\right)
                         _{i,j\in\{0,1,\dots,k-1\}}
  \end{equation}
  is the number of weighted families of non-intersecting lattice paths from
  $$(0,0),(2,-2),\dots,\left(2(k-1),-2(k-1)\right)$$ 
  to
  $$(3,m-k-1),(5,m-k-2),\dots,\left(2(k-1)+3,m-2k\right),$$
  where a vertical step with an even $x$-coordinate has weight $q$ and all
  other steps have weight $1$.
\end{thm}
\begin{proof}
  The determinantal formula follows from the preceding lemma. The combinatorial
  interpretation is a standard application of the main theorem of
  nonintersecting lattice paths, and completely analogous to the applications
  given in Gessel and Viennot{\cite{GesselViennot1989}}.
\end{proof}
\begin{cor}
  The coefficients of $P_{m,k}$ are nonnegative and symmetric.
\end{cor}
\begin{proof}
  A combinatorial way to see the symmetry is as follows: Modifying the weights
  such that vertical steps with an odd $x$-coordinate have weight $q$ and all
  the others weight $1$ does not change the entries of the determinant.
  
  However, consider any given family of paths with weight $q^w$, when vertical
  steps with even $x$-coordinate have weight $q$. After the modification of the
  weights it will have weight $q^{max-w}$, where $max$ is the total number of
  vertical steps in such a family of paths, which implies the claim.
\end{proof}
\begin{rmk}
  It appears that the polynomials $P_{m,k}$ are log-concave, however, we did
  not pursue this question further.
\end{rmk}

\section{Acknowledgements}

Many thanks are due to Michael Schlosser and Christian Krattenthaler for their
patience and help with Lemma~\ref{lem:1} and for pointing out some bad typos in
the manuscript.  
%\bibliography{math} 
\bibliographystyle{amsplain}

\end{document}